\documentclass[a4paper]{amsart}

\usepackage{amsfonts, latexsym, amssymb, graphicx,
multicol, mathrsfs, color, amscd, verbatim, paralist,
xspace, url, euscript, stmaryrd,  amsmath, amscd, enumitem,
bbold, multirow, tikz, mathtools, mathabx}
\usepackage[all,pdf,cmtip]{xy}

\usepackage{chngcntr}
\usepackage{apptools}

\usepackage{xpatch}
\xapptocmd{\appendix}{%
  \counterwithin{theorem}{section}
  \counterwithin{definition}{section}
  \counterwithin{table}{section}
}{\typeout{Success}}{}

\usepackage[colorlinks, linkcolor=blue, citecolor=magenta, urlcolor=cyan]{hyperref}

\usepackage{mathabx}

\def\qed{{\hfill $\Box$}}
\newtheorem{theorem}{THEOREM}[section]
\newtheorem{corollary}[theorem]{Corollary}

\theoremstyle{definition}

\theoremstyle{remark}
\newtheorem{remark}[theorem]{Remark}

\newcommand\CC{{\mathbb C}}
\newcommand\RR{{\mathbb R}}

\def\Re{\mathop{\rm Re}\nolimits}
\def\Im{\mathop{\rm Im}\nolimits}

\makeatletter
\def\blfootnote{\xdef\@thefnmark{}\@footnotetext}
\makeatother

\begin{document}

\title[Real-analyticity of rigid spherical hypersurfaces]{On the real-analyticity of rigid
\vspace{0.1cm}\\
spherical hypersurfaces in $\CC^2$}

\author[Isaev]{Alexander Isaev}
\author[Merker]{Jo\"el Merker}

\address{Mathematical Sciences Institute\\
Australian National University\\
Canberra, Acton, ACT 2601, Australia}
\email{alexander.isaev@anu.edu.au}

\address{Laboratoire de Math\'ematiques d'Orsay\\
Universit\'e Paris-Sud\\
CNRS, Universit\'e Paris-Saclay\\
91405 Orsay Cedex\\
France}
\email{joel.merker@math.u-psud.fr}

\maketitle

\thispagestyle{empty}

\pagestyle{myheadings}

\begin{abstract}
We prove that every smooth rigid spherical hypersurface in $\CC^2$ is in fact real-analytic. As an application of this result, it follows that the classification of real-analytic rigid spherical hypersurfaces in $\CC^2$ found by V.~Ezhov and G.~Schmalz applies in the smooth case.
\end{abstract}

\section{Introduction}\label{intro}
\setcounter{equation}{0}

We consider connected smooth real hypersurfaces in the space $\CC^2$ with standard coordinates $z=x+iy,w=u+iv$. For simplicity, we understand smoothness as $C^{\infty}$-smoothness although everything we do works for the class $C^k$ with sufficiently large $k$, say for the class $C^7$. Specifically, we look at \emph{rigid hypersurfaces}, i.e., hypersurfaces given by equations of the form
\begin{equation}
u=F(z,\overline{z}),\label{rigideqgen}
\end{equation}
where $F$ is a smooth real-valued function defined on a domain $U\subset\CC$. Rigid hypersurfaces are invariant under the 1-parameter family of holomorphic transformations
\begin{equation}
(z,w)\mapsto (z,w+it),\,\, t\in\RR.\label{translatisymm}
\end{equation}
Throughout the paper we assume Levi-nondegeneracy, which for a rigid hypersurface  means that $F_{z\bar z}$ is everywhere nonvanishing (here and below $z,\bar z$ used as subscripts indicate the corresponding partial derivatives). 


All our considerations will be entirely local, in a neighborhood of a chosen point. In particular, letting 0 be the reference point, we utilize a natural notion of local equivalence for germs of rigid hypersurfaces. Namely, two rigid hypersurface germs at 0 are called \emph{rigidly equivalent} if there exists a map of the form
$$
\begin{array}{l}
(z,w)\mapsto (g(z),a w+h(z) ),\,\,a\in \RR^*
\end{array}
$$
that transforms one germ into the other, where $g, h$ are holomorphic near 0 and $g'(0)\ne 0$, $g(0)=h(0)=0$.

In this paper, we focus on the following problem: 
\vspace{-0.7cm}\\

$$
\begin{array}{l}
\hspace{0.2cm}\hbox{$(*)$ classify, up to rigid equivalence, the germs of rigid hypersurfaces that are}\\
\vspace{-0.3cm}\hspace{0.7cm}\hbox{\emph{spherical}, i.e., CR-equivalent to a germ of the sphere $S^3\subset\CC^2$.}\\
\end{array}
$$
\vspace{-0.4cm}\\

The sphericity condition for a Levi-nondegenerate smooth hypersurface $M\subset\CC^3$ is equivalent to that of the vanishing of the CR-curvatu\-re, which is a ${\mathfrak{su}}(2,1)$-valued 2-form defined on an 8-dimensional principal fiber bundle over $M$. This form arises from the reduction of the CR-structures of 3-dimensional Levi-nondegenerate CR-manifolds of hypersurface type to absolute parallelisms performed in \cite{Ca} and later generalized to higher dimensions, in particular, in papers \cite{T1}, \cite{T2}, \cite{T3}, \cite{CM}, \cite{Ch}, \cite{BS2}, \cite{BS3}. The interpretation of the sphericity condition as that of CR-flatness is a basis of our approach to task $(*)$. 

The problem of describing CR-flat structures that possess certain symmetries (e.g., as in (\ref{translatisymm})) is a natural one and has been addressed by a number of authors. In particular, under the assumption of CR-flatness, \emph{homogeneous} strongly pseudoconvex CR-hypersurfaces have been studied (see \cite{BS1}), and \emph{tube} Levi-nondegenerate hypersurfaces in complex space have been extensively investigated and even fully classified for certain signatures of the Levi form (see \cite{I2} for a detailed exposition). Furthermore, in \cite{I3} a class of CR-flat \emph{Levi-degenerate} tube hypersurfaces in $\CC^3$ was fully described. Compared to the tube case, the case of rigid hypersurfaces is the next situation up in terms of complexity and is substantially harder. For instance, as explained in \cite{I4}, the study of Levi-degenerate rigid hypersurfaces in $\CC^3$ is a much more difficult task in comparison with the work on tube hypersurfaces done in \cite{I3}. 

Although not mentioned explicitly, problem $(*)$ was first looked at in article \cite{S} for \emph{real-analytic} hypersurfaces. Even in this simplified setup, determining all rigid spherical hypersurfaces turned out to be highly nontrivial, with only a number of examples found in \cite{S}. A complete solution to $(*)$ in the real-analytic category was only recently obtained in \cite{ES} by employing techniques based on the rigid normal forms constructed in \cite{S}. Although the hypersurface germs on the list found in \cite{ES} are given by implicit equations, these equations are not very complicated and represent an acceptable solution to problem $(*)$ in the real-analytic situation. For the completeness of our exposition, we state the result of \cite{ES} in Section \ref{esclass} (see Theorem \ref{resultes}).

In this note we show that the classification of \cite{ES} applies in the smooth case as well. Namely, we obtain the following:

\begin{theorem}\label{main}
Every smooth rigid spherical hypersurface in $\CC^2$ is real-analytic.
\end{theorem}

Our proof of Theorem \ref{main} in the next section is based on the zero CR-curvature equations for rigid Levi-nondegenerate hypersurfaces in $\CC^n$ for any $n\ge 2$ found in \cite{I1} (see \cite[Theorem 2.1]{I2} for a better exposition). These equations form a system of PDEs for the graphing function, which turns into a single equation for $F$ in $\CC^2$ (see formula (\ref{eqspher}) in Theorem \ref{thmrem} below). We note that another PDE for the function $F$, also characterizing the sphericity property of a rigid Levi-nondegenerate hypersurface in $\CC^2$, was produced in \cite{L} (see formula (\ref{alternativeeq})). The argument of \cite{L} is based on the Chern-Moser normal forms, so it assumes real-analyticity. In fact, equation (\ref{alternativeeq}) can be shown to characterize sphericity even in the smooth setup, hence it is equivalent to (\ref{eqspher}). One way to see this is by specializing calculations in articles \cite{MS1}, \cite{MS2} to the rigid case. Thus, each of (\ref{eqspher}), (\ref{alternativeeq}) can be used to establish Theorem \ref{main}, but the argument based on equation (\ref{eqspher}) is shorter, so we chose to utilize (\ref{eqspher}) in our proof. We sketch a proof relying on equation (\ref{alternativeeq}) in Remark \ref{alternativeproof}. 

{\bf Acknowledgements.} This work was initiated while the first author was visiting D\'epartement de Math\'ematiques d'Orsay and completed while he was visiting the Steklov Mathematical Institute in Moscow. We thank S. Nemirovski and A. Domrin for useful discussions.

\section{Proof of Theorem \ref{main}}\label{proof}
\setcounter{equation}{0}

We start by stating the special case of \cite[Theorem 2.1]{I2} in complex dimension two. 

\begin{theorem}\label{thmrem} Consider a rigid hypersurface $M$ in $\CC^2$ given by equation {\rm (\ref{rigideqgen})} with the function $F$ defined on a domain $U\subset\CC$. Assume that $F$ is everywhere Levi-nondegenerate, i.e., $F_{z\bar z}\ne 0$ on $U$. Then $M$ is spherical if and only if $F$ satisfies a differential equation of the form
\begin{equation}
F_{zz}=A(F_z)^3+B(F_z)^2+CF_z+D,\label{eqspher}
\end{equation}
where $A$, $B$, $C$, $D$ are functions holomorphic on $U$.
\end{theorem}

Set $f:=F_z$. To prove Theorem \ref{main} it suffices to show that $f$ is real-analytic on $U$. Indeed, if $f$ were real-analytic, then locally near every $z_0\in U$ we would have
$$
F=\int f dz+\varphi(\bar z),\label{exprforf}
$$
where $\int f dz$ denotes the result of the term-by-term integration with respect to $z-z_0$ of the power series in $z-z_0$, $\overline{z-z_0}$ representing $f$  near $z_0$, and $\varphi$ is antiholomorphic hence real-analytic. It would then follow that $F$ is real-analytic at $z_0$ and therefore on $U$.

To prove that $f$ is real-analytic on $U$, observe that (\ref{eqspher}) implies        
\begin{equation}
f_z=Af^3+Bf^2+Cf+D.\label{eqsphernew}
\end{equation}
Separating the real and imaginary parts of equation (\ref{eqsphernew}), we see that it is equivalent to a system of two equations for two real-valued functions:
\begin{equation}
\begin{array}{l}
r_x+s_y=G(z,\bar z,r,s),\\
\vspace{-0.3cm}\\
\hspace{-0.3cm}-r_y+s_x=H(z,\bar z,r,s),
\end{array}\label{ellipticsys}
\end{equation}
where $r:=\Re f$, $s:=\Im f$ and $G$, $H$ are real-valued analytic functions of all their arguments. Clearly, system (\ref{ellipticsys}) is nonlinear in general. Since
$$
\det\left(
\begin{array}{rl}
\lambda_1 & \lambda_2\\
-\lambda_2 & \lambda_1
\end{array}
\right)=\lambda_1^2+\lambda_2^2
$$
is positive for all nonzero vectors $(\lambda_1,\lambda_2)\in\RR^2$, system (\ref{ellipticsys}) is elliptic as defined in \cite[pp.~210, 266]{M2}. Then by \cite[p.~203]{M1} (see also \cite[Theorem 6.7.6]{M2}), every solution of (\ref{ellipticsys}) is real-analytic. We thus see that $f$ is real-analytic as required.\qed  
 
\begin{remark}\label{alternativeproof}
As mentioned in the introduction, in article \cite{L} another equation characterizing the sphericity property of a real-analytic rigid Levi-nondegenerate hypersurface in $\CC^2$ was obtained, and, by results of \cite{MS1}, \cite{MS2}, this equation works in the smooth setup as well. Let
$$
\mu:=\log(F_{z\bar z})_{\bar z}=\frac{F_{z\bar z\bar z}}{F_{z\bar z}}.
$$
Then the equation found in \cite{L} is as follows:    
\begin{equation}
\mu_{z\bar z\bar z}-3\mu_{z\bar z}\mu+2\mu_z\mu^2-\mu_z\mu_{\bar z}=0.\label{alternativeeq}
\end{equation}

We will now show that (\ref{alternativeeq}) leads to an elliptic system, just as (\ref{eqsphernew}) does. Indeed, separating the real and imaginary parts of equation (\ref{alternativeeq}), we see that it is equivalent to a system of two equations for two real-valued functions:
\begin{equation}
\begin{array}{l}
\nu_{xxx}+\nu_{xyy}-\eta_{xxy}-\eta_{yyy}=\Phi,\\
\vspace{-0.3cm}\\
\nu_{xxy}+\nu_{yyy}+\eta_{xxx}-\eta_{xyy}=\Psi,
\end{array}\label{ellipticsys1}
\end{equation}
where $\nu:=\Re \mu$, $\eta:=\Im \mu$ and $\Phi$, $\Psi$ are analytic functions of $\nu$, $\eta$, as well as their first and second partial derivatives with respect to $x$, $y$. To see that system (\ref{ellipticsys1}) is elliptic, we compute
$$
\det\left(
\begin{array}{lr}
\lambda_1^3+\lambda_1\lambda_2^2 & -\lambda_1^2\lambda_2-\lambda_2^3\\
\vspace{-0.3cm}\\
\lambda_1^2\lambda_2+\lambda_2^3 & \lambda_1^3-\lambda_1\lambda_2^2
\end{array}
\right)=\lambda_1^6+\lambda_1^4\lambda_2^2+\lambda_1^2\lambda_2^4+\lambda_2^6,
$$
which is positive for all nonzero vectors $(\lambda_1,\lambda_2)\in\RR^2$. Hence, as before, by \cite[p.~203]{M1} and \cite[Theorem 6.7.6]{M2} every solution of (\ref{ellipticsys1}) is real-analytic, and therefore $\mu$ is real-analytic. It then follows that locally near every $z_0\in U$ we have
$$
F_{z\bar z}=\exp\left(\int\mu d\bar z\right)\exp\left(\varphi(z)\right),
$$
where $\varphi$ is holomorphic. Thus, $F_{z\bar z}$ is real-analytic at $z_0$, and, by the ellipticity of Poisson's equation, $F$ is real-analytic on $U$. 
\end{remark}

\section{The Ezhov-Schmalz classification}\label{esclass}
\setcounter{equation}{0}

In this section, for the completeness of our exposition, we state the main result of \cite{ES}, which is a solution to problem $(*)$ in the real-analytic setup.

\begin{theorem}\label{resultes}
Every germ at the origin of a real-analytic rigid spherical hypersurface in $\CC^2$ is rigidly equivalent to the germ at the origin of the rigid hypersurface defined by an equation of the form
\begin{equation}
\begin{array}{l}
\displaystyle(1+2\phi|z|^2)\frac{\sin(2ru)}{2r}-e^{-2\theta u}|z|^2-\\
\vspace{-0.5cm}\\
\hspace{1.5cm}\displaystyle\left(\phi+\bar c z+c \bar z+4\phi(\phi-\theta)|z|^2\right)\frac{e^{-2\theta u}-\cos(2ru)+\frac{\theta\sin(2ru)}{r}}{r^2+\theta^2}=0.
\end{array}\label{esclassform}
\end{equation}
Here $\phi,\theta\in\RR$, $c\in\CC$, $r\in\RR\cup i\RR$ and there exists a pair of numbers $\tau,\rho\in\RR$ such that:
\begin{itemize}

\item [] $\phi$ is any real root of the cubic polynomial
\begin{equation}
4\phi^3+4\tau\phi^2+(\tau^2-\rho)\phi-|c|^2,\label{eqforphi}
\end{equation}
and

\item[]
$$
\begin{array}{l}
\theta=\tau+3\phi,\\
\vspace{-0.1cm}\\
r^2=3\phi^2+2\tau\phi-\rho.
\end{array}
$$

\end{itemize} 

\noindent In formula {\rm (\ref{esclassform})}, by expanding into power series and taking limits, one sets
$$
\displaystyle\frac{\sin(2ru)}{2r}\Bigg|_{r=0}:=u
$$ 
and
$$
\frac{e^{-2\theta u}-\cos(2ru)+\frac{\theta\sin(2ru)}{r}}{r^2+\theta^2}\Bigg|_{r=\theta=0}:=2u^2.
$$
\end{theorem}

Notice that every equation in (\ref{esclassform}), up to multiplying $r$ by -1 and up to the three possible choices of $\phi$ as a real root of (\ref{eqforphi}), is uniquely determined by three numbers: $c\in\CC$, $\tau\in\RR$, $\rho\in\RR$. Multiplication of $r$ by -1 leads to a rigidly equivalent hypersurface germ, whereas it is not clear whether different choices of $\phi$ can result in some nonequivalent germs. However, even if polynomial (\ref{eqforphi}) has a unique real root, two distinct triplets $(c,\tau,\rho)$ and $(c',\tau',\rho')$ do not necessarily yield rigidly nonequivalent hypersurface germs. For example, for $c=\tau=0$ and any $\rho<0$ we see $\phi=\theta=0$, $r=\pm\sqrt{-\rho}$. In this case formula (\ref{esclassform}) becomes
$$
\frac{\sin(\pm2\sqrt{-\rho}\,u)}{\pm2\sqrt{-\rho}}-|z|^2=0. 
$$
By scaling $z,w$ one observes that the germ at 0 of this hypersurface for any $\rho<0$ is rigidly equivalent to the germ of the hypersurface  
$$
\sin u=|z|^2,
$$
which is one of the examples found in \cite{S}. Thus, Theorem \ref{resultes} does not describe the moduli space of real-analytic rigid spherical hypersurface germs in $\CC^2$ precisely.
\vspace{0.1cm}

We conclude our paper with the following consequence of Theorem \ref{main}:

\begin{corollary}\label{corol}
Every germ at the origin of a smooth rigid spherical hypersurface in $\CC^2$ is rigidly equivalent to the germ at the origin of the rigid hypersurface defined by an equation of the form {\rm (\ref{esclassform})} as described in the statement of Theorem {\rm\ref{resultes}}.
\end{corollary}


\begin{thebibliography}{ABC}

\bibitem[BS1]{BS1} Burns, D. and Shnider, S., Spherical hypersurfaces in complex manifolds, {\it Invent. Math.} {\bf 33} (1976), 223--246.

\bibitem[BS2]{BS2} Burns, D. and Shnider, S., Real hypersurfaces in complex manifolds, in {\it Several Complex Variables, Proc. Sympos. Pure Math.} {\bf XXX}, Amer. Math. Soc., Providence, R.I., 1977, pp. 141--168.

\bibitem[BS3]{BS3} Burns, D. and Shnider, S., Projective connections in CR geometry, {\it Manuscripta Math.} {\bf 33} (1980/81), 1--26.

\bibitem[Ca]{Ca} Cartan, \' E., Sur la g\'eometrie pseudo-conforme des hypersurfaces de l'espace de deux variables complexes: I, {\it Ann. Math. Pura Appl.} {\bf 11} (1933), 17--90; II, {\it Ann. Scuola Norm. Sup. Pisa} {\bf 1} (1932), 333--354.

\bibitem[CM]{CM} Chern, S. S. and Moser, J. K., Real hypersurfaces in complex manifolds, {\it Acta Math.} {\bf 133} (1974), 219--271; erratum, {\it Acta Math.} {\bf 150} (1983), 297.

\bibitem[Ch]{Ch} Chern, S. S., On the projective structure of a real hypersurface in $\CC^{n+1}$, {\it Math. Scand.} {\bf 36} (1975), 74--82.

\bibitem[ES]{ES} Ezhov, V. and Schmalz, G., Explicit description of spherical rigid hypersurfaces in $\CC^2$, {\it Complex Anal. Syner.} {\bf 1}:2 (2015), DOI: 10.1186/2197-120X-1-2.

\bibitem[I1]{I1} Isaev, A. V., Rigid spherical hypersurfaces, {\it Complex Var. Theory Appl.} {\bf 31} (1996), 141--163.

\bibitem[I2]{I2} Isaev, A. V., {\it Spherical Tube Hypersurfaces}, Lect. Notes Math. {\bf 2020}, Springer, New York, 2011.

\bibitem[I3]{I3} Isaev, A. V., Affine rigidity of Levi degenerate tube hypersurfaces, {\it J. Differential Geom.} {\bf 104} (2016), 111--141.

\bibitem[I4]{I4} Isaev, A. V., Rigid Levi degenerate hypersurfaces with vanishing CR-curvature, {\it J. Math. Anal. Appl.} {\bf 474} (2019), 782--792.

\bibitem[L]{L} Loboda, A. V., On the sphericity of rigid hypersurfaces in $\CC^2$, {\it Math. Notes} {\bf 62} (1997), 329--338.

\bibitem[MS1]{MS1} Merker, J. and Sabzevari, M., Explicit expression of Cartan's connection for Levi-nondegenerate 3-manifolds in complex surfaces, and identification of the Heisenberg sphere, {\it Central Europ. J. Math.} {\bf 10} (2012), 1801--1835.

\bibitem[MS2]{MS2} Merker, J. and Sabzevari, M., The Cartan equivalence problem for Levi-nondegenerate real hypersurfaces $M^5\subset\CC^2$, {\it Izv. Math.} {\bf 78} (2014), 1158--1194.

\bibitem[M1]{M1} Morrey, C. B., On the analyticity of the solutions of analytic non-linear elliptic systems of partial differential equations: Part I, analyticity in the interior, {\it Amer. J. Math.} {\bf 80} (1958), 198--218.

\bibitem[M2]{M2} Morrey, C. B., {\it Multiple Integrals in the Calculus of Variations}, Springer, Berlin, 2008.

\bibitem[S]{S} Stanton, N., A normal form for rigid hypersurfaces in $\CC^2$, {\it Amer. J. Math.} {\bf 113} (1991), 877--910.

\bibitem[T1]{T1} Tanaka, N., On generalized graded Lie algebras and geometric structures I, {\it J. Math. Soc. Japan} {\bf 19} (1967), 215--254.

\bibitem[T2]{T2} Tanaka, N., On non-degenerate real hypersurfaces, graded Lie algebras and Cartan connections, {\it Japan. J. Math.} {\bf 2} (1976), 131--190.

\bibitem[T3]{T3} Tanaka, N., On the equivalence problems associated with simple graded Lie algebras, {\it Hokkaido Math. J.} {\bf 8} (1979), 23--84.

\end{thebibliography}
\end{document}